\documentclass[12pt]{amsart}%
\usepackage{graphicx}
\usepackage{psfrag}
\usepackage{mathrsfs}
\usepackage{amsmath,amssymb}
\font\de=cmssi12

\begin{document}


\newtheorem{conj}{\sc Conjecture}
\newtheorem{defi}{\sc Definition}

\newtheorem{cor}{\sc Corrollary}
\newtheorem{ques}{\sc Question}[section]
\newtheorem{teo}{\sc Theorem}
\renewcommand{\theteo}{\Alph{teo}}
\newtheorem{teor}[teo]{\sc Theorem }
\newtheorem{prop}{\bf \sc Proposition}[teo]
\newtheorem{propa}{\bf \sc Proposition}
\renewcommand{\thepropa}{1.\arabic{propa}}
\newtheorem{propb}{\bf \sc Proposition}
\renewcommand{\thepropb}{2.\arabic{propb}}
\newtheorem{obs}{\sc Remark}[section]
\newtheorem{lema}{\sc Lemma}[prop]
\newtheorem{lemaa}{\sc Lemma}
\renewcommand{\thelemaa}{1.\arabic{lemaa}}
\newtheorem{lemab}{\sc Lemma}
\renewcommand{\thelemab}{2.\arabic{lemab}}
\newtheorem{claim}{\it Claim}[prop]
\renewcommand{\theprop}{\Alph{teo}.\arabic{prop}}
\def\bp{\noindent{\it Proof. }}
\def\ep{\noindent{\hfill $\fbox{\,}$}\medskip\newline}
\renewcommand{\theequation}{\arabic{section}.\arabic{equation}}
\newcommand{\thetheorem}{\arabic{section}.\arabic{theorem}}
\newcommand{\eps}{\varepsilon}
\newcommand{\disp}[1]{\displaystyle{\mathstrut#1}}
\newcommand{\fra}[2]{\displaystyle\frac{\mathstrut#1}{\mathstrut#2}}
\newcommand{\dif}{{\rm Diff}^r}
\newcommand{\ph}{{\rm PH}^r_m}
\newcommand{\phr}[1]{{\rm PH}^{\mathstrut#1}_m}
\newcommand{\Per}{{\rm Per}}
\newcommand{\Z}{\mathbb Z}
\newcommand{\R}{\mathbb R}
\newcommand{\N}{\mathbb N}
\newcommand{\s}{\sigma}
\newcommand{\W}{\mathcal W}
\newcommand{\F}{\mathcal F}
\newcommand{\D}{\mathscr D}
\newcommand{\KS}{\mathscr{K\!\!S}}
\newcommand{\stack}[2]{\tiny{\hspace*{-.5em}\begin{array}{l}#1\\\noalign{\vspace*{-.3em}}#2\end{array}}}
\def\to{\mathop{\rightarrow}}
\def\ord{\mathop{\rm ord}}
\def\diam{\mathop{\rm diam}}
\def\cc{\mathop{\rm cc}}
\title[Accessibility and stable ergodicity]{Accessibility and stable ergodicity for partially hyperbolic diffeomorphisms with 1D-center bundle}
\author{F.Rodriguez Hertz}
\author{M.Rodriguez Hertz}
\author{R. Ures}
\thanks{This work was partially supported by FCE 9021, CONICYT-PDT 29/220 and  CONICYT-PDT 54/18 grants}
\address{IMERL-Facultad de Ingenier\'\i a\\ Universidad de la
Rep\'ublica\\ CC 30 Montevideo, Uruguay} \email{frhertz@fing.edu.uy}
\email{jana@fing.edu.uy} \email{ures@fing.edu.uy}
\maketitle
\begin{abstract}We prove that stable ergodicity is $C^r$ open and dense among
conservative partially hyperbolic diffeomorphisms with
one-dimensional center bundle, for all $r\in[2,\infty]$.\par The
proof follows Pugh-Shub program \cite{pughshub2000}: among
conservative partially hyperbolic diffeomorphisms with
one-dimensional center bundle, accessibility is $C^r$ open and
dense, and essential accessibility implies
ergodicity.\end{abstract}
\section{Introduction}

In the second half of the 19th century Boltzmann introduced the
term {\em ergodic} within the context of the study of gas
particles and since then, even in its initial formulation the
Ergodic Hypothesis was extremely unlikely, ergodic theory grew up
to be a useful tool in many branches of physics.\par Subsequent
reformulations and developments turned the original ergodic
hypothesis into the statement: {\em time average equals space
average} for typical orbits, that is
$$\lim_{n\to\infty}\frac1n\sum_{k=0}^{n-1}\phi(f^k(x))=\int_M\phi\,d\mu\qquad \mu-{\rm a.e.} x$$
A system is {\de $\mu$-ergodic} if it satisfies the hypothesis
above for all $C^0$ observables $\phi$, or equivalently, if only
full or null $\mu$-volume sets are invariant under the dynamics.
Near 1930, after the first ergodic theorems appeared -
\cite{neumann}, \cite{birkhoff}, \cite{birkkoop} - it was
conjectured that most conservative systems were ergodic.
\par With the Kolmogorov-Arnold-Moser (KAM) phenomenon (1954) it came
out that there were full open sets of conservative non-ergodic
systems \cite{kam}. Indeed, KAM theory presented completely
integrable systems, a dynamic that could be described as elliptic,
for which a big amount (positive volume) of invariant tori survived
after performing perturbations, which prevents ergodicity. This
is an example of a {\em stably non ergodic} system.\par %
On the other end of the spectrum, the work of Hopf \cite{hopf}, and
later Anosov-Sinai \cite{anosov,anosovsinai}, gave full open sets of
ergodic systems, a fact that was unknown up to that time. Anosov
systems, are what we call now completely hyperbolic dynamics, and
were for some time the only {\de stably ergodic} examples known. By
stably ergodic is meant a diffeomorphism in the interior of the set
of ergodic diffeomorphisms.
\par
Almost three decades later, Grayson, Pugh, Shub got the first
non-hyperbolic example of a stably ergodic system \cite{gps}.
These examples have a {\em partially hyperbolic dynamics}
\cite{bp}, \cite{hps}: there are strong contracting and strong
expanding invariant directions, but a center direction also
appears. Since then, the area became quite active and many stably
ergodic examples appeared, see \cite{pughshub2004} for a survey.
Let us also mention that there are already examples of
conservative stably ergodic systems that are {\em not} partially
 hyperbolic \cite{ali}. \par%
In this new context, Pugh and Shub have proposed the
following:\vspace*{.6em}\newline
{\sc Conjecture 1 } {\em Stable ergodicity is $C^r$ open and dense
among volume preserving partially hyperbolic diffeomorphisms, for
all $r\geq 2$.} \vspace{.5em}\par
As far as we know, the conjecture above was first stated in 1995,
at the International conference on dynamical systems held in
Montevideo, Uruguay \cite{pughshub1995}. We thank Keith Burns for
this information. \par
In this paper, we prove this conjecture is true in case the center
bundle is one dimensional: \setcounter{teo}{-1}
\begin{teo}[\sc Main]\label{stable ergodicity is open and dense}
Stable ergodicity is $C^r$ open and dense among volume preserving
partially hyperbolic diffeomorphisms with one dimensional center
distribution, for all $r\geq 2$.
\end{teo}
 In \cite{pughshub2000}, Pugh and Shub proposed a program for the
 proof of this conjecture. This approach was based on the notion
 of {\em accessibility}: A diffeomorphism $f$ has the
 {\de accessibility property} if the only non void set consisting of
 whole
  stable leaves and whole unstable leaves is the  manifold $M$ itself. It
   has the {\de essential accessibility property} if every measurable
set consisting of whole stable leaves and  whole unstable leaves has
 full or null volume. Clearly, accessibility implies essential
accessibility. When talking about stable and unstable leaves we are
referring to the leaves of the unique foliations tangent to the
contracting and expanding directions, respectively.\par%
Pugh and Shub suggested the following two conjectures:
\vspace*{.6em}\newline
{\sc Conjecture 2:} {\em Stable accessibility property is open and
dense among $C^r$ partially hyperbolic diffeomorphisms, volume
preserving or not, $r\geq 2$.} \vspace{.8em}\par
In the case $\dim E^c=1$, the accessibility property is always
stable \cite{didier}. For the sake of simplicity, let us call
$\ph(M)$ the set of partially hyperbolic $C^r$ diffeomorphisms of
$M$, preserving a smooth probability measure $m$. In this paper, we
prove that:
\begin{teo}\label{accessibility is open and dense}
Accessibility is open and dense in $\ph(M)$, for all $1\leq r\leq
\infty$, if the center distribution is one dimensional.
\end{teo}
In fact, we obtain that accessibility is $C^1$ open and $C^\infty$
dense. Let us observe that the conjecture is established here only
for the conservative case. Earlier results in this direction can
be found in \cite{niticatorok}, where they prove stable
accessibility is $C^r$ dense for one-dimensional center bundle,
under certain hypotheses (for instance, dynamical coherence and
compact center leaves), and \cite{dolgowilk}, where stable
accessibility is shown to be dense in the $C^1$ topology with no
assumption on the dimension of the center bundle.
\par The second conjecture of the Pugh-Shub program is:\vspace*{.5em}
\newline
{\sc Conjecture 3:} {\em Essential accessibility implies ergodicity
among $C^2$ volume preserving partially hyperbolic diffeomorphisms.}
\vspace{.5em}\newline
We also prove this conjecture in case the center dimension is one.
\begin{teo}\label{essential accessibility implies kolmogorov}
Essential accessibility implies Kolmogorov (in particular,
ergodicity) in $\phr2(M)$, if the center distribution is
one-dimensional.
\end{teo}\vspace{-.1em}
Let us mention that K. Burns and A. Wilkinson have recently proved a
result that implies theorem \ref{essential accessibility implies
kolmogorov}: they show that essential accessibility implies
Kolmogorov in $\phr2(M)$ under the assumption of a mild center
bunching condition, with no assumption on the dimension of the
center bundle. No dynamical coherence is required \cite{bw2}. They
also prove that differentiability condition in Theorem
\ref{essential accessibility implies kolmogorov} can be improved to
$C^{1+\mbox{\tiny H\"{o}lder}}$. We thank A. Wilkinson
for this information.\par%
{\bf Acknowledgements.} We want to thank M. Shub for his support in
a difficult moment. We also want to thank K. Burns for reading early
versions of this manuscript and for useful remarks. We are also
grateful to C. Pugh for many valuable suggestions.

\section{Preliminaries, notation and sketch of the
proof}\label{section preliminaries and notation}
Let $M$ be a compact Riemannian manifold, and $m$ be a smooth
probability measure on $M$. Denote by $\dif_m(M)$ the set of $C^r$
volume preserving diffeomorphisms. In what follows we shall
consider a {\de partially hyperbolic} $f\in\dif_m(M)$, that is, a
diffeomorphism admitting a non trivial $Df$-invariant splitting of
the tangent bundle $TM=E^s\oplus E^c\oplus E^u$, such that all
unit vectors $v^\s\in E^\s_x$ ($\s=s,c,u$) with $x\in M$ verify:
$$\|T_xfv^s\|<\|T_xfv^c\|<\|T_xfv^u\|$$
for some suitable Riemannian metric. It is also required that
$\|Tf|_{E^s}\|<1$ and $\|Tf^{-1}|_{E^u}\|<1$. We shall denote by
$\ph(M)$ the family of $C^r$ volume preserving partially hyperbolic
diffeomorphisms of $M$.
\par It is a known fact that there are foliations ${\mathcal
W}^\s$ tangent to the distributions $E^\s$ for $\s=s,u$ (see for
instance \cite{bp}). A set $X$ will be called {\de $\s$-saturated}
if it is a union of leaves of ${\mathcal W}^\s$, $\s=s,u$.\par%
In this paper we will consider the case $\dim E^c=1$. After Peano,
we can find small curves passing through each $x\in M$, that are
everywhere tangent to the bundle $E^c$. We shall call these curves
{\de center curves through $x$}, and denote them by Roman
$W^c_{loc}(x)$, since a priori they are non-uniquely integrable
curves, in order to distinguish them from the true foliations
$\W^\s$, $\s=s,u$. It is easy to see that
$f$ takes center curves into center curves.\par %
We shall denote by $\W^\s(x)$ the leaf of $\W^\s$ through $x$
for ($\s=s,u$) and will write $\W^\s_{loc}(x)$ for a small disk in
$\W^\s(x)$ centered in $x$. For any choice of $W^c_{loc}(x)$, the
sets
$$W^{\s c}_{loc}(x)=\W^\s_{loc}(W^c_{loc}(x))=\bigcup_{y\in W^c_{loc}(x)}\W^\s_{loc}(y)\qquad \s=s,u$$
are $C^1$ (local) manifolds everywhere tangent to the sub-bundles
$E^\s\oplus E^c$ for $\s=s,u$ (see, for instance
\cite{brin.burago.ivanov}). The sets above depend on the choice of
$W^c_{loc}(x)$.
\begin{obs}\label{remark center en cs}
Observe that for all choices of $W^{sc}_{loc}(x)$ and $y\in
W^{sc}_{loc}(x)$, there exists a center curve $W^c_{loc}(y)$
through $y$ contained in $W^{sc}_{loc}(x)$ (see
\cite{brin.burago.ivanov})
\end{obs}
\subsection{Proof of theorem \ref{accessibility is open and dense}}
Let us say that a set $\Gamma$ is {\de $\s$-saturated} if $\Gamma$
is union of leaves of ${\mathcal W}^\s$, $\s=s,u$. For the proof of
theorem \ref{accessibility is open and dense}, we will see that
$C^r$-generically, the minimal $s$- and $u$- saturated set that
contains any point $x$ (that is, the {\de accessibility class} of
$x$) is the whole $M$. This property is known as the {\de
accessibility property} and is open in $\phr1(M)$ if the center
bundle is one-dimensional \cite{didier}.\par The proof focuses on
the open accessibility classes, and the first step is showing that
for any periodic point, a perturbation can be made so that its
accessibility class becomes open (Unweaving Lemma). Secondly, we
obtain periodic points for any dynamics in $\ph(M)$ having non
trivial open accessibility classes that do not cover $M$. A
genericity argument allows us to conclude, via Kupka-Smale
techniques the following statement:
\setcounter{teo}{1}
\begin{prop}\label{proposition set C(M)}
$C^r$-generically in $\ph(M)$, $r\geq2$, either one of the
 following properties holds:
\begin{enumerate} \item $f$ has the accessibility property or \item\label{su integrable sin pper} $
\Per(f)=\emptyset$ and the distribution $E^s\oplus E^u$ is
integrable
\end{enumerate}
\end{prop}
As one would expect the second possibility is quite unstable under
perturbations and, indeed, this is the case:
\begin{prop}\label{set B is nowhere dense}
Situation {\rm(\ref{su integrable sin pper})} described above is
nowhere dense in $\ph(M)$.
\end{prop}
We show that the Unweaving Lemma mentioned above holds also for non
recurrent points. In this way, integrability of $E^s\oplus E^u$ can
be broken by small perturbations.\par%
In both cases, to have some control on how perturbations affect
local invariant manifolds, we need the existence of points whose
orbits keep away from the support of the perturbation (Keepaway
Lemma \ref{keepaway}).
\par The two statements
together imply theorem \ref{accessibility is open and dense}. This
part is developed in \S\ref{section accessibility is open and
dense}.
\subsection{Proof of Theorem \ref{essential accessibility implies
kolmogorov}}%
\label{subsection proof of teorema essential accessibility}
For the proof of Theorem \ref{essential accessibility implies
kolmogorov}, we shall mainly follow the line in \cite{gps},
\cite{pughshub2000} and \cite{bw}. This theorem was obtained
independently of \cite{bw2}, though Burns and Wilkinson's result is
more general. We decided to include Theorem \ref{essential
accessibility implies kolmogorov} here for completeness, and because
it is simpler in the sense that it uses true leaves instead of fake
foliations, which are a difficult (and possibly necessary if $\dim
E^c>1$) technical step. Also, it takes two steps to characterize
Lebesgue density points instead of the seven equivalences in \S4 of
\cite{bw2}. \begin{ques} Is it possible to use the techniques here
and avoid the fake foliations in case the bunching conditions in
\cite{bw2} hold and $E^c$ is weakly integrable, that is there are
center leaves everywhere tangent to $E^c$ at every point?
\end{ques} Let us consider a diffeomorphism $f$
having the {\de essential accessibility} property, that is,
verifying that each measurable $s$- and $u$-saturated set is full or
null measure. In order to prove that $f$ is {\de ergodic} (each {\em
invariant} set is full or null measure) it suffices to show, due to
Birkhoff's ergodic theorem, that
$$\phi_{\pm}(x)=\limsup_{n\to\infty}\frac1n\sum_{k=1}^n\phi(f^{\pm
k}(x))=\int_M\phi\,{\rm dm}\qquad m \,\, {\rm  a.e.} x$$ for all
$C^0$ observables $\phi:M\to \R$. It is not hard to see that, for
each $c\in \R$, the set $S(c)=\phi_+^{-1}(c,\infty)$ is
$s$-saturated, and the set $U(c)=\phi_-^{-1}(c,\infty)$ is
$u$-saturated. Since $m(S(c)\triangle U(c))=0$ due to Birkhoff's
theorem, we have that the set $S(c)\cap U(c)$ differs in a set of
null measure from an $s$- saturated set, and also from a
$u$-saturated set. In general, we shall say that a measurable set
$X$ is {\de essentially $\s$-saturated} if there exists a measurable
$\s$-saturated set $X_\s$ (an {\de essential $\s$-saturate of $X$})
such that $m(X\triangle X_\s)=0$. In short, $S(c)\cap U(c)$ is
essentially $s$- and essentially $u$-saturated (with essential
$s$-and $u$- saturates $S(c)$ and $U(c)$, respectively).\par The
typical Hopf's argument went on by showing that in fact Lebesgue
density points of any set $X$ were $s$- and $u$-saturated, whence
the essential accessibility property directly implied ergodicity.
The differentiability of holonomy maps played an important role in
this fact. However, in this context we {\em do not have}
differentiable holonomy maps.\par This gap will be covered by
proving instead that \setcounter{teo}{2}\setcounter{prop}{0}
\begin{prop}\label{prop2.2}
The Lebesgue density points of any essentially $s$- and
essentially $u$-saturated set $X$ form an $s$- and $u$-saturated
set.
\end{prop}
That is, Lebesgue density points of essentially $s-$ and essentially
$u-$saturated sets flow through stable and unstable leaves. In
\cite{pughshub2000}, Pugh and Shub suggested that certain shapes
called {\de juliennes} would be more natural, rather than merely
Riemannian balls, in order to treat preservation of density points.
Here we follow this line and use certain solid juliennes instead of
balls.\par%
Of course, these new neighborhood bases will define different sets
of density points. We will consider the following generalization of
Lebesgue density points:\par Let us say that a point $x$ is a {\de
$C_n$-density point} of a set $X$ if $\{C_n(x)\}_n$ is a local
neighborhood basis of $x$, and
$$\lim_{n\to \infty}\frac{m(X\cap C_n(x))}{m(C_n(x))}=1$$
In particular, the {\de Lebesgue density points} will be the
$\{B_{r^n}(x)\}_{n\geq1}$-density points, where $B_{r^n}(x)$ is the
Riemannian ball centered at $x$ with radius $r^n$, $r\in(0,1)$. The
choice of $r$ is irrelevant, since $x$ is a $B_{r^n}$-density point
of $X$ if and only if $$\lim_{\eps\to 0}\frac{m(X\cap
B_\eps(x))}{m(B_\eps(x))}=1.$$ A {\de $cu$-julienne} $J^{cu}_n(x)$
of $x$ is a dynamically defined local unstable saturation of a
center curve, its radius depending on $x$ and $n$, and going to $0$
subject to certain rates related to contraction rates in the bundles
(see precise definitions in \S\ref{subsection definition of
juliennes}, formulas (\ref{definition juliennes})). We shall define
a {\de solid julienne} $J^{suc}_n(x)$ of $x$ as a local stable
sat\-u\-ra\-tion of some $cu$-julienne (precise definitions in
\S\ref{section cu density points}). The family
$\{J^{suc}_n(x)\}_{n\geq1}$ is a measurable neighborhood basis of
$x$. For this family we obtain
\begin{prop}\label{saturacion de puntos de densidad}
The set of $J^{suc}_n$-density points of an essentially
$s$-saturated set $X$ is $s$-saturated.
\end{prop}
By changing the neighborhood basis, we have solved the problem of
preserving density points, that is we have established Proposition
\ref{prop2.2} but for julienne density points. However, we need to
know now what the relationship is between the julienne density
points, and Lebesgue density points. Given a family $\mathcal{M}$ of
measurable sets, let us say that two systems $\{C_n\}_n$ and
$\{E_n\}_n$ are {\de Vitali equivalent} over $\mathcal{M}$, if the
set of $C_n$-density points of $X$ equals (pointwise) the set of
$E_n$-density points of $X$ for all $X\in\mathcal{M}$. The argument
is completed by showing that
\begin{prop}\label{cu density points}
The family $\left\{J^{suc}_n(x)\right\}$ is Vitali equivalent to
Lebesgue over essentially $u$-saturated sets, for any choice of
$W^c_{loc}(x)$.
\end{prop}

Hence, over essentially $s$- and $u$-saturated sets, the set of
Lebesgue density point is $s$- saturated. A symmetric argument
shows it is also $u$-saturated.\par%
This ends the proof of Proposition \ref{prop2.2} and, actually, it
shows essential accessibility implies ergodicity. To show that, in
fact, it implies Kolmogorov property, \cite{pesin} states that it
suffices to see that the Pinsker algebra (the largest subalgebra for
which the entropy is zero) is trivial. But after \cite{bp}, sets in
the Pinsker algebra are essentially $s$- and essentially
$u$-saturated, what proves Theorem \ref{essential accessibility
implies kolmogorov}.

\section{Accessibility is $C^r$ open and dense}\label{section
accessibility is open and dense}
Let us call $AC(x)$ the accessibility class of the point $x$. We
will show that the set
$$\D=\left\{f\in\ph(M): AC(x) \mbox{ is open for all $x\in \Per(f)$}\right\}$$
is $C^r$ dense, where $\Per(f)$ denotes the set of periodic points
of $f$. This is the set $\D$ mentioned in Proposition
\ref{proposition set C(M)}. Afterwards, as stated in that
proposition, it will be shown that $\D$ may be decomposed into a
disjoint union
\begin{equation}\label{C=AUB}
\D={\mathscr A}\cup{\mathscr B}\end{equation} where ${\mathscr
A}(M)$ consists of diffeomorphisms with the accessibility property
and ${\mathscr B}$ consists of diffeomorphisms lacking  periodic
points and verifying that the distribution $E^s\oplus E^u$ is
integrable. Moreover, ${\mathscr B}$ will be shown to be nowhere
dense. This will prove Proposition \ref{set B is nowhere dense} and,
in fact, Theorem \ref{accessibility is open
and dense}.\par%
In this section, we shall denote, for any set $X\subset M$,
$$\W^\s_{loc}(X)=\bigcup_{x\in X}\W^\s_{loc}(x)\qquad
\mbox{with}\quad\s=s,u.$$
\subsection{A lamination in the complement of open accessibility classes}
Fix $f\in \ph(M)$, and let $U(f)$ be the set of points whose
accessibility classes are open, and $\Gamma(f)=M\setminus U(f)$ be
the complement of $U(f)$. A {\de lamination} ${\mathcal L}$ is a
foliation of a closed subset $N\subset M$. In this case, we say $N$
is {\de laminated} by ${\mathcal L}$.
\setcounter{teo}{1}\setcounter{prop}{2}
\begin{prop}\label{lamina}
$\Gamma(f)$ is a compact, invariant set laminated by the
accessibility classes.
\end{prop}
\begin{prop}\label{abierto}
For a given point $x\in M$ the following statements are equivalent
\begin{enumerate}
\item $AC(x)$ is open.%
\item $AC(x)$ has non empty interior.%
\item $AC(x)\cap W^c_{loc}(x)$ has nonempty interior for any
choice of $W^c_{loc}(x)$.%
\end{enumerate}
\end{prop}
\noindent%
Indeed, an open set within the accessibility class may be joined
to any other point $z$ of $AC(x)$ by an {\de $su$-path}: a path
consisting of a finite number of arcs, each contained either in an
$s$- or a $u$-leaf (see figure).
\begin{figure}[h]
\includegraphics[width=12cm]{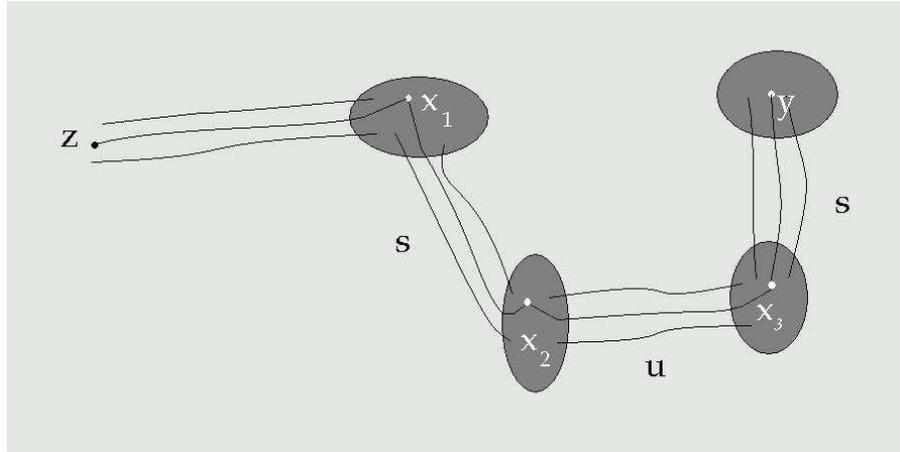}
\caption{An $su$ path from $z$ to $y$}
\end{figure}
Let $z=x_0, x_1,\dots,x_{n-1}, x_n=y$ be points in the $su$-path
such that $x_i$ and $x_{i+1}$ are in the same $\s$-leaf (for $\s$
either $s$ or $u$). And let $y$ be in the interior of $AC(x)$. The
$\s$- holonomy maps are continuous, so there exists a neighborhood
of $x_n$ contained in $AC(x)$. A finite inductive argument allows us
to conclude $z$ is also in the interior of $AC(x)$, so $AC(x)$ is
open.\par
\begin{figure}[h]
\includegraphics[width=9cm]{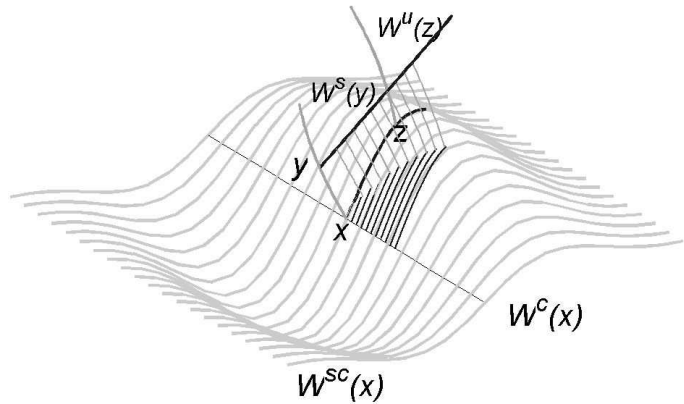}
\caption{\label{acopen proyeccion}A point in $U(f)$ (open
accessibility class)}
\end{figure}

%
If $AC(x)$ is open, it is obvious that $AC(x)\cap W^c_{loc}(x)$
will have non empty interior for any choice of $W^c_{loc}(x)$, but
the converse is also true. Indeed, there is a well defined map
$$p_{us}:W^{usc}_{loc}(x)\to W^c_{loc}(x)$$ where $W^{usc}_{loc}(x)=\W^u_{loc}(W^{sc}_{loc}(x))$, that is obtained by
first projecting along $\W^u$ and then along $\W^s$ (see Figure
\ref{acopen proyeccion}). If a point $w$ is in an open set $V$ of
$AC(x)\cap W^c_{loc}(x)$, then $p_{us}^{-1}(V)$ will be an open
neighborhood of $w$, due to continuity of $p_{us}$. But
$p_{us}^{-1}(V)$ is clearly in $AC(x)$, hence $AC(x)$ is open.\ep%
Let $AC_x(y)$ denote the connected component of $AC(y)$ in
$W^{usc}_{loc}(x)$ containing $y$. The points in $\Gamma(f)$ have
the following property:
\begin{lema}
Let $z\in\Gamma(f)$, then for all points $y\in AC(z)\cap
W^{usc}_{loc}(x)$ and all $w\in p_{us}^{-1}(x)$, the set $AC_x(y)$
meets $W^c_{loc}(w)$ at exactly one point.
\end{lema}
\begin{proof}
$p_{us}(AC_x(y))$ is a connected subset of $W^c_{loc}(x)$. Thus
$W^c_{loc}(x)$, being one dimensional, would have non empty interior
if it contained more than one point, which would contradict
Proposition \ref{abierto}. As $p_{us}$ restricted to any
$W^c_{loc}(x)$ is one to one, we get the proposition.
\end{proof}
\begin{figure}[b]
\vspace{-.5cm}
\includegraphics[width=10cm]{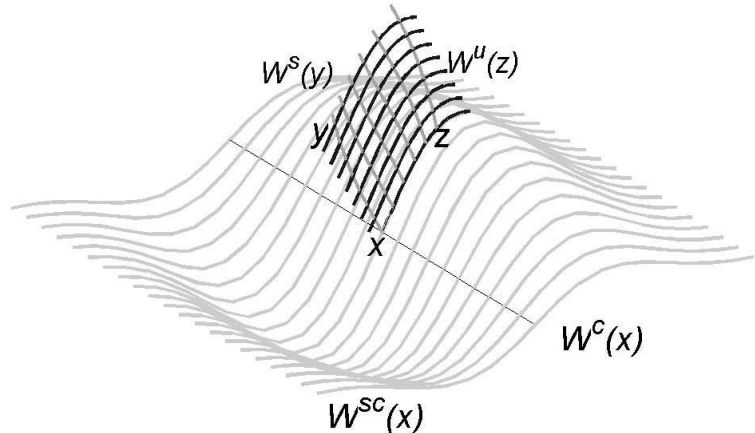}
\vspace{-.5cm} \caption{A point in $\Gamma(f)$}
\end{figure} Finally, from the preceding lemma we get that, if we denote
$\Gamma^c_x=\Gamma(f)\cap W^c_{loc}(x)$, then the set
$\W^u_{loc}(\W^s_{loc}(\Gamma^c_x))$ is exactly $\Gamma(f)\cap
W^{usc}_{loc}(x)$, and for any point $y\in\Gamma^c_x$, we have
$AC_x(y)=\W^u_{loc}(\W^s_{loc}(y))$. It is not hard to see that the
lamination charts are coherent, and hence we get Proposition
\ref{lamina}.%
\begin{obs} \label{remark su integrable}
Observe that, in particular, if $AC(x)$ is not open then
$$\W^u_{loc}(y)\cap\W^s_{loc}(z)\ne\emptyset$$
for all $y\in\W^s_{loc}(x)$ and all $z\in\W^{u}_{loc}(x)$.
\end{obs}

\subsection{Keepaway Lemma}

Let  $f$ be a diffeomorphism  preserving a foliation $\W$ tangent to
a continuous sub-bundle $E\subset TM$.  Call $\W(x)$ the leave of
$\W$ through $x$ and $\W_\eps(x)$ the set of points  that are
reached from $x$ by a curve contained in $\W(x)$ of length less than
$\eps$. If $V=V(x)$ is a (small) disk trough $x$ transverse to $\W$
whose dimension equals  the codimension of $E$, define
$B_\eps(V)=\cup \{\W_\eps(y);\; y\in V\}$ and
$C_\eps(V)=B_{4\eps}(V)\setminus B_\eps(V)$.

\begin{lema}[Keepaway Lemma]\label{keepaway} Suppose that under the previous conditions we have, in
addition, that $||Tf^{-1}|_E||<\mu^{-1}<1$. Let $N$ be such that
$\mu^N>4$. If there exist $x \in M$, $V(x)$ and $\varepsilon
>0$  such that:
$$f^n(C_\eps(V))\cap B_{\varepsilon}(V)=\emptyset\qquad \forall n=1,\dots ,N$$
 then for every $y\in V$ there
exists $z\in \W_{4\varepsilon}(y)$ such that $f^n(z)\notin
B_{\varepsilon}(V)\quad \forall n>0$.
\end{lema}

\bp Let $y\in V$ and $w\in \W_{4\eps}(y)$ such that
$\W_\eps(w)\subset C_\eps(V)$. Call $D_0=\overline{\W_\eps(w)}$. We
shall construct, by induction, a sequence of closed disks $D_n$ such
that $f^{-1}(D_n)\subset D_{n-1}\;\forall n>0$ and $D_n \cap
B_\eps(V)=\emptyset$. Thus $z$ will be any point in $\cap
\{f^{-n}(D_n);\; n\in \mathbb{N}\}$ (in fact in our construction
this intersection will consist in a unique point).

The construction is as follows:

\begin{enumerate}
\item \label{primerpaso} If $n< N$ then $D_n=f^n(D_0)$.
\item \label{segundopaso} There exists a point $w_N\in f^N(D_0)$ such that $\W_{4\eps}(w_N)\subset
f^N(D_0)$. Let $D_N=\overline{\W_{4\eps}(w_N)}$. Observe that, for
$n>N$, $\W_{4\eps}(f^{n-N}(w_N))\subset f^{n-N}(D_N)$.
\item Suppose that $n>N$ is such that $\W_{4\eps}(f^{j-N}(w_N))\cap B_\eps(V)=\emptyset
$ for $j=N,\dots, n$. Then define
$D_n=\overline{\W_{4\eps}(f^{n-N}(w_N))}$.
\item If $n_1$ is such that $\W_{4\eps}(f^{j-N}(w_N))\cap B_\eps(V)=\emptyset
$ for $j=N,\dots, n_1-1$ and $\W_{4\eps}(f^{n_1-N}(w_N))\cap
B_\eps(V)\neq\emptyset $ we have that there exists a point
$w_{n_1}\in \W_{4\eps}(f^{n_1-N}(w_N))$ such that
$\overline{\W_\eps(w_{n_1})}\subset C_\eps (V)$. Define
$D_{n_1}=\overline{\W_\eps(w_{n_1})}$.
\item Now, to continue the construction, go to step
\ref{primerpaso}, and substitute $D_0$ by $D_{n_1}$.
\end{enumerate}

This algorithm gives the desired sequence of disks, and then the
point $z$, proving the lemma. \ep

\begin{obs}\label{obsnonreturning}
Sometimes the following consequence is more useful than Lemma
\ref{keepaway}:

Let $x$, $V(x)$ and $\eps$ be as in Lemma \ref{keepaway}.  Let $y$
and $\delta>0$ be such that $f^i(\W_\delta(y))\cap
B_\eps(V(x))=\emptyset$ for $i=1,\dots,K$ where $K$ is such that
$\mu^K\delta>4\eps$. Then, there is $z\in \W_\delta(y)$ so that
$f^n(z)\notin B_{\varepsilon}(V)\quad \forall n>0$.
\end{obs}

\bp Observe that $\W_{4\eps}(f^K(y))\subset f^K(\W_\delta(y))$. Now
go to step \ref{segundopaso} in the algorithm of the lemma replacing
$w_N$ by $f^K(y)$ . \ep

Call $\mathscr I=\{f\in \ph{M};\, E^s\oplus E^u \mbox{ is
integrable} \}$. Observe that  $\mathscr I$ is a closed set and
$\mathscr B\subset \mathscr I$.\par
In the partially hyperbolic setting the Keepaway Lemma
\ref{keepaway} and Remark \ref{obsnonreturning} have as corollaries
that $\mathscr I$ has empty interior and that, given a periodic
point $x$, $f$ can be perturbed in such a way that the accessibility
class of $x$ for the perturbed diffeomorphism is open. This is shown
in the next subsections.
\subsection{$\D$ is dense}
Genericity of $\D$ follows from the classical Kupka-Smale argument,
after the following property:

\begin{lema}[Unweaving Lemma] \label{unweaving} For each $x\in \Per(f)$ there exists $g$ $C^r$-close
to $f$ such that $x\in \Per(g)$ and $AC_g(x)$ is open.
\end{lema}
\bp Assume that $AC_f(x)$ is not open for some periodic point
$x\in\Per(f)$. Then, as stated in remark \ref{remark su integrable},
$\W^s_{loc}(y)\cap\W^u_{loc}(z)\ne \emptyset$ for all
$y\in\W^u_{loc}(x)$ and all $z\in \W^s_{loc}(x)$.\par The idea is to
perturb a small neighborhood of $x$, so that $x\in\Per(g)$ and
$\W^s_{g,loc}(\hat y)\cap \W^u_{g, loc}(\hat z)= \emptyset$ for some
$\hat y\in \W^u_{g,loc}(x)$ and $\hat z\in \W^s_{g,loc}(x)$. This
will obviously prove $AC_g(x)$ is open.\par

Observe that, taking $V(x)=W^{sc}_{loc}(x)$  for some $W^c_{loc}(x)$
and $\eps>0$ small, the point $x$ verifies the hypothesis of the
Keepaway Lemma \ref{keepaway}. Then, we obtain a point $y\in
\W_{4\eps}^u(x)$ such that its forward orbit does not intersect
$B_\eps(V(x))$. Analogously, applying the Keepaway Lemma to
$f^{-1}$, we obtain a point $z\in \W_{4\eps}^s(x)$ that does not
return for the past to a similar neighborhood of $x$, say $B_\eps
(\hat V(x))$. Now, we can choose $k>0$ and a small $\delta
>0$ in such a way that $\W_\delta^s(f^{-k}(y))$,
$\W_\delta^u(f^{k}(y))$  and $\{w\}=\W_\delta^u(f^{k}(z))\cap
\W_\delta^s(f^{-k}(y))$ are contained in $B_\eps=B_\eps(V(x))\cap
B_\eps (\hat V(x))$. Call $\hat y=f^{-k}(y)$ and $\hat z=f^{-k}(z)$.
\par

From the way in which $y$ and $z$ are chosen we can take $U$, a
sufficiently small  neighborhood of $w$, in such a way that $f^n
(\W_\delta^s(\hat y))$ and $f^{-n} (\W_\delta^u(\hat z))$ does not
cut $U$ for all $n>0$. Also we can require $U$ not to intersect
$\W^\sigma_{\eps}(f^n(x))$ for all $n$, $\sigma=u,s$.

  It follows that $\W^s_{\eps}(x)$,
$\W^u_{\eps}(x)$, $\W^u_\delta(\hat z)$ and $f(\W_\delta^s(\hat y))$
do not change if we perform a perturbation supported in $U$. Now it
is easy to perturb $f$ in $U$ so that $g(\W_\delta^u(\hat z)\cap
U)\cap f(\W_\delta^s(\hat y))=\emptyset$.
This implies that $\W^s_{g,loc}(\hat y)\cap \W^u_{g, loc}(\hat z)= \emptyset$ and finishes the proof of the lemma.%
\ep
\begin{figure}[b]
\hspace*{-.5cm}\includegraphics[width=8cm]{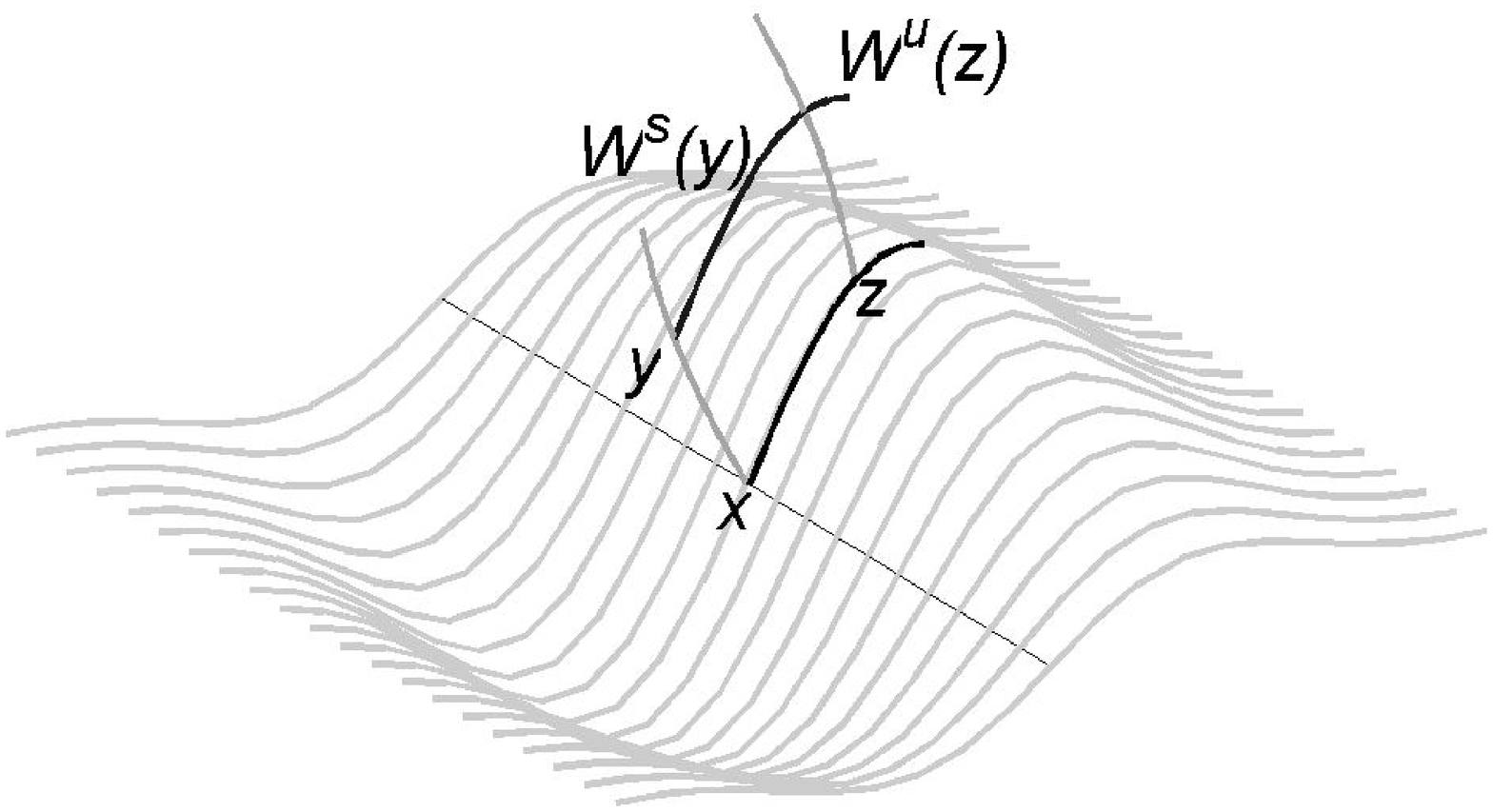}\includegraphics[width=8cm]{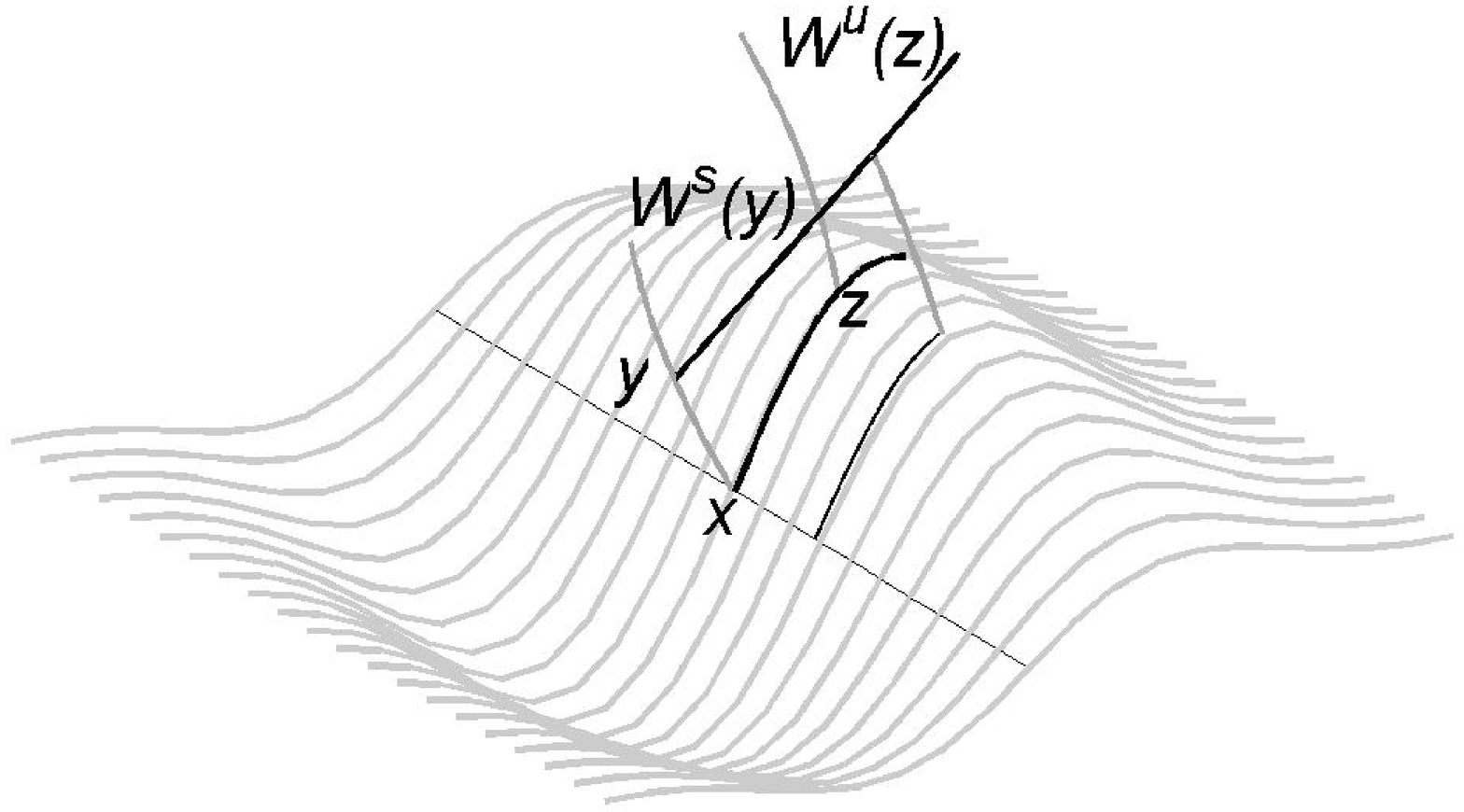}
\caption{Unweaving Lema: Before and after perturbing around a
periodic point $x$}
\end{figure}

The Unweaving Lemma above implies, after Kupka-Smale, that
$C^r$-generically it holds:
$$\Per(f)\subset U(f)$$
This means, the set $\D$ is $C^r$-generic. The following proposition
shows that, in case $\Gamma(f)$ is a proper subset, there are always
periodic points in $\Gamma(f)$. This situation is nowhere dense.
\begin{prop}\label{gama tiene periodicos}
If $\emptyset\varsubsetneq\Gamma(f)\varsubsetneq M$, then
$Per(f)\cap \Gamma(f)\neq\emptyset$.
\end{prop}
\bp Let us prove there is a periodic point in the boundary $\partial
\Gamma(f)$ of $\Gamma(f)$. Observe that $\partial \Gamma(f)$ is a
compact, $f$-invariant, $su$-saturated set. We will assume $M$ and
$E^c$ are orientable. Indeed, by taking a double covering if
necessary, we can assume $M$ is orientable. If $E^c$ is not
orientable, we take again a double covering $\tilde M$ of $M$ in
such a way that $\tilde E^c$, the lift of $E^c$, is orientable. Let
$\tilde f$ be a lift of $f$ to $\tilde M$, then ${\tilde f}^2$ is
partially hyperbolic, $\tilde E^c$ is its center bundle and $\tilde
f^2$ preserves the orientation of $\tilde E^c$. Any point $x\in
\Gamma(f)$ lifts to a point ${\tilde x}\in {\tilde \Gamma}(\tilde
f^2)\subset {\tilde M}$. The set ${\tilde \Gamma}(\tilde f^2)$ is
locally diffeomorphic to
$\Gamma(f)$, and is ${\tilde f}^2$ invariant. So we shall assume that $M$ and $E^c$ are orientable.\par%

Take a point $x\in \partial\Gamma(f)$. We may also assume, without
loss of generality, that there is an open interval
$I=(x,x+\triangle x)^c$ contained in $W^c_{loc}(x)\setminus
\Gamma(f)$ with $x+\triangle x\notin\Gamma(f)$. Let us call
$V=\W^u_{loc}(\W^s_{loc}(I))$, so $f^k(V)\cap \Gamma(f)=\emptyset$
for all $k\in\Z$. Observe that, if we denote by $(a_y,b_y)^c$ the
component of $y$ in the set $W^c_{loc}(y)\cap V$, then $a_y$ is
always in $AC_x(x)\subset\partial\Gamma(f)$, and $b_y$ is never in
$\Gamma(f)$.\par%
Now, as the non-wandering set of $f$ is $M$, there exists $y\in V$
such that $f^k(y)\in V$ for some $k>0$. Indeed, $f^k(a_y)\in
\partial \Gamma(f)$, then Lemma \ref{lema integrable} together with
the fact that $f$ preserves the orientation of $E^c$ imply that
$f^k(a_y)$ actually belongs to $AC_x(x)$ (see also Figure \ref{non.uniquely.integrable}).\par%
The proof follows now from the standard lemma:
\begin{lema}\label{lema creacion periodicos}
There is $\eps_0>0$ such that if $x\in\Gamma(f)$ verifies
$f^k(B^{su}_{\eps_0}(x))\cap B^{su}_{\eps_0}(x)\neq\emptyset$ for
some $k>0$, then there is a periodic point in
$B^{su}_{\eps_0}(x)$.
\end{lema}\ep
The following property is a consequence of continuity and
transversality of the invariant bundles, and has been used in
proving Proposition \ref{gama tiene periodicos}:
\begin{figure}[h]
\includegraphics[height=6cm]{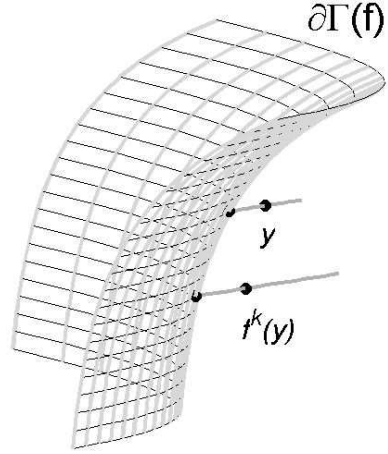}
\caption{\label{non.uniquely.integrable}Bounded dispersion of
center curves}
\end{figure}
\begin{lema}\label{lema integrable}
For each small $\eps>0$ there exists $\delta>0$ such that if
$d(x,y)<\delta$ and $z\in W^c_\delta (x)$, then $W^c_{loc}(y)\cap
\W^s_\eps(\W^u_\eps(z))\ne\emptyset$, regardless of the choice of
center leaves of $x$ and $z$.
\end{lema}
After Proposition \ref{gama tiene periodicos}, we have the following
possibilities for $f\in\D$:
\begin{enumerate}
\item $\Gamma(f)=\emptyset$, that is, $f$ has the accessibility
property%
\item $\Gamma(f)=M$ with $\Per(f)=\emptyset$
\end{enumerate}
The situation $\emptyset\varsubsetneq \Gamma(f)\varsubsetneq M$
cannot happen for $f\in \D$, since it implies there is a
periodic point in $\Gamma(f)$. This proves Proposition \ref{proposition set C(M)}%

\subsection{Proposition \ref{accessibility is open and dense}.\ref{su integrable sin pper}}

Recall that $\mathscr I=\{f\in \ph(M);\, E^s\oplus E^u \mbox{ is
integrable} \}$ and that since $\mathscr I$ is a closed set and
$\mathscr B\subset \mathscr I$, next proposition implies Proposition
\ref{accessibility is open and dense}.\ref{su integrable sin pper}.

\begin{prop}$\mathscr I(M)$ has empty interior.
\end{prop}

\bp The proof is similar to that of the Unweaving Lemma
\ref{unweaving}. Take any nonperiodic point $x\in M$. Given $N>0$
there exists $\varepsilon>0$ such that $f^i (B_\eps(V(x)))\cap
B_{\varepsilon}(V(x))=\emptyset$ $\forall i=0,\dots ,N,\, i\neq 0$.
Then,  Remark \ref{obsnonreturning} applied to $x$ itself implies
that there exists $z\in \W^u_{\varepsilon}(x)$ such that
$f^n(z)\notin B_{\varepsilon}(V(x))\, \forall n> 0$. Take $\gamma>0$
such that $\W^s_{\gamma}(z)\subset B_{\varepsilon}(x)$ and
$f^{n+1}(\W^s_{\gamma}(z))\cap B_{\varepsilon}(x)= \emptyset \,
\forall n\geq 0$.

By applying Remark \ref{obsnonreturning} three times we obtain a
nonrecurrent point $w\in W^s_{\gamma}(z)$, points
$\overline{x},\overline{y}$ (different from $w$) and $\rho>0$ such
that:

\begin{itemize}
\item $\overline{x}\in \W^u_{\rho}(w)$ and $f^n(\overline{x})\notin
B_{\rho}(w) \, \forall n>0$

\item $\overline{y}\in \W^s_{\rho}(w)$ and
$f^{-n}(\overline{y})\notin B_{\rho}(w) \, \forall n>0$
\end{itemize}

Let $\overline{z}=\W^s_{loc}(\overline{x})\cap
\W^u_{loc}(\overline{y})$.

Finally observe that if we perform a perturbation in a small
neighborhood of $f^{-1}(\overline{z})$ we have that
$\W^s_{loc}(w),\, \W^u_{loc}(w),\, \W^s_{loc}(\overline{x})$  and
$\W^u_{loc}({f^{-1}(\overline{y})})$ remain unchanged. Similarly to
the proof of the Unweaving Lemma \ref{unweaving}, we can do this
perturbation in order to obtain that
$\W^{u}_{g,\,loc}(\overline{y})$ does not intersect $\W^{s}_{g,\,
loc}(\overline{x})$ for the perturbed diffeomorphism $g$. \ep

\section{Essential accessibility implies ergodicity}\label{section
ess acc implies erg}
\subsection{Definitions}\label{subsection definition of juliennes}
 Let us consider smooth functions $\nu, \hat\nu,\gamma,\hat\gamma:M\to \R^+$
 verifying, for all unit vectors $v^i\in E^i$ with $i=s,c,u$ and $x\in
 M$,
$$\|T_xfv^s\|<\nu(x)<\gamma(x)<\|T_xfv^c\|<\hat\gamma(x)^{-1}<\hat\nu(x)^{-1}<\|T_xfv^u\|$$
where $\nu,\hat\nu<1$ and $\|.\|$ is an adapted Riemannian metric as
at the beginning of the section. We may also assume that $d$ and
$\nu,\hat\nu, \gamma, \hat\gamma$ verify:
\begin{equation}\label{adapted.metric.su}
\begin{array}{ll}
 d(f(x),f(x'))\leq\nu(x)\,d(x,x')&\mbox{ for}\quad
x'\in\W^s_{loc}(x)\\\noalign{\medskip} d(f^{-1}(x),f^{-1}(x'))\leq
\hat\nu(f^{-1}(x))\,d(x,x')&\mbox{ for}\quad x'\in\W^u_{loc}(x)
\end{array}
\end{equation}
\begin{equation}\label{adapted.metric.c}
\begin{array}{ll}
d(f(x),f(x'))\leq \hat\gamma(x)^{-1}d(x,x')&\mbox{ for}\quad x'\in
W^c_{loc}(x)\\\noalign{\medskip}
d(f^{-1}(x),f^{-1}(x'))\leq\gamma(f^{-1}(x))^{-1}d(x,x')& \mbox{
for}\quad x'\in W^c_{loc}(x)
\end{array}
\end{equation}
with
\begin{equation}\label{bunching condition}
\begin{array}{c}
\fra{\nu(x)}{\gamma(x)}<\sigma<\min (1,{\hat\gamma(x)})
 \end{array}
\end{equation} for some
smooth $\sigma:M\to\R$. Note that $\nu,\hat\nu<1$, while
$\gamma.\hat\gamma$ and $\sigma$ can be chosen less than but close
to $1$.
\begin{obs}
Inequalities (\ref{adapted.metric.su}) and (\ref{adapted.metric.c})
do not depend on the choice of the center curve through $x$.
\end{obs}
Consider, for $\alpha=\nu, \hat\nu,\gamma, \hat\gamma,\s$ and $n\geq
0$ the multiplicative cocycles:
$$\alpha_n(x):=\prod_{i=0}^{n-1}\alpha(f^i(x))\qquad \alpha_{-n}(x):=\alpha_n(f^{-n}(x))^{-1}$$
For each $W^c_{loc}(x)$, define the set
$$B^c_n(x)=W^c_{\s_n(x)}(x)$$ and consider also:
\begin{equation}\label{definition juliennes}
J^u_n(x)=f^{-n}(\W^u_{\nu_n(x)}(f^n(x)))\qquad\mbox{and}\qquad
J^{cu}_n(x)=\bigcup_{y\in B^c_n(x)}J^u_n(y)
\end{equation}
The sets $J^{cu}_n(x)$ will be called center-unstable juliennes of
$x$ or {\de $cu$-juliennes}
\subsection{Controlling stable holonomy}\label{subsection stable
holonomy}
In this section we will prove that the deformation suffered by the
$cu$-juliennes under the stable holonomy, can be controlled in a the
following sense:\setcounter{teo}{2}
\begin{prop}\label{holonomia controlada}
There exists $k\in{\mathbb Z}^+$ such that, if $x'\in
\W^s_{loc}(x)$, then for all choices of $W^c_{loc}(x)$ and
$W^c_{loc}(x')$ contained in $W^{sc}_{loc}(x)$, the stable
holonomy map from $W^{cu}_{loc}(x)$ to $W^{cu}_{loc}(x')$ verifies
$$J^{cu}_{n+k}(x')\subset h^s(J^{cu}_n(x))\subset J^{cu}_{n-k}(x')\qquad \forall \,n\geq k$$
\end{prop}
The proof splits into two parts. On one hand, we prove that the
holonomy does not distort center leaves too much, as it is seen in
Lemma \ref{holonomia controlada sobre center} and Figure \ref{figura
holonomia central}. On the other hand, it is seen that each unstable
fiber on a certain center leaf, is transformed, under the stable
holonomy in a curve contained in a greater julienne. This is seen in
Lemma \ref{holonomia no se tuerce sobre inestable} and Figure
\ref{figura inestable no se tuerce}.
\begin{figure}[h]
\includegraphics[width=11cm]{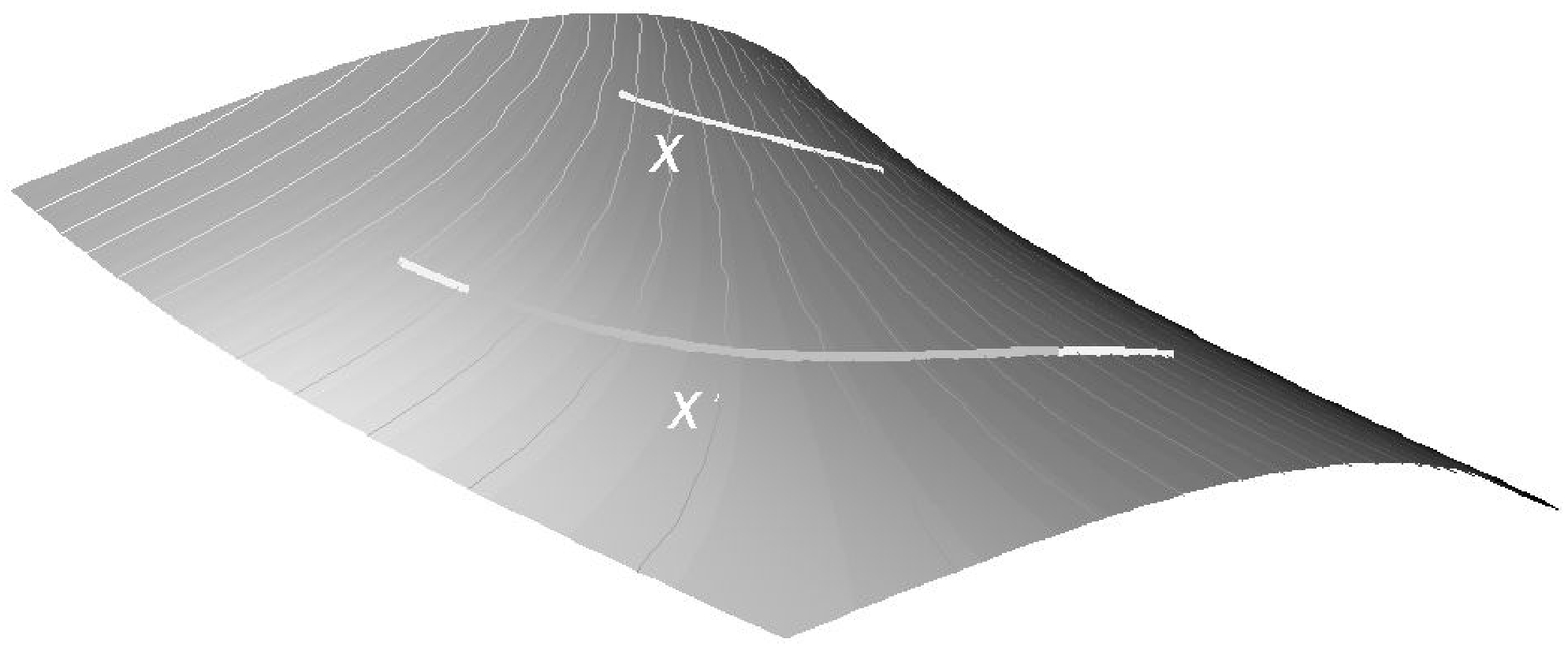}
\caption{\label{figura holonomia central}Lemma \ref{holonomia
controlada sobre center}}
\end{figure}
\begin{lema}\label{holonomia controlada sobre center}
There exists $k\in\Z^+$, not depending on $x$, such that for all
choices $W^c_{loc}(x)$, $W^c_{loc}(x')$ of center curves through
$x$, $x'$ contained in some $W^{sc}_{loc}(x)$, with
$x'\in\W^s_{loc}(x)$, the stable holonomy map $h^s$ from
$W^c_{loc}(x)$ to $W^c_{loc}(x')$, verifies
$$h^s(B^c_n(x))\subset B^c_{n-k}(x')\qquad \forall n\geq k$$
\end{lema}
\bp Consider $L>0$  and $C>1$ be as in Proposition \ref{holonomia
derivable} and Lemma \ref{holder cocycles} of Appendix 1,
respectively. Take $k>0$ such that $\s_{-k}(x)>LC$ for all $x\in M$
(recall that $\s<1$), then
$$h^s(B_n^c(x))\subset W^c_{L\s_n(x)}(x')\subset W^c_{LC\s_n(x')}(x')\subset W^c_{\s_{n-k}(x')}(x')=B^c_{n-k}(x')$$
and the claim follows.\ep%
The following lemma is the second part of the proof of proposition
\ref{holonomia controlada}:
\begin{figure}[h]
\includegraphics[width=11cm]{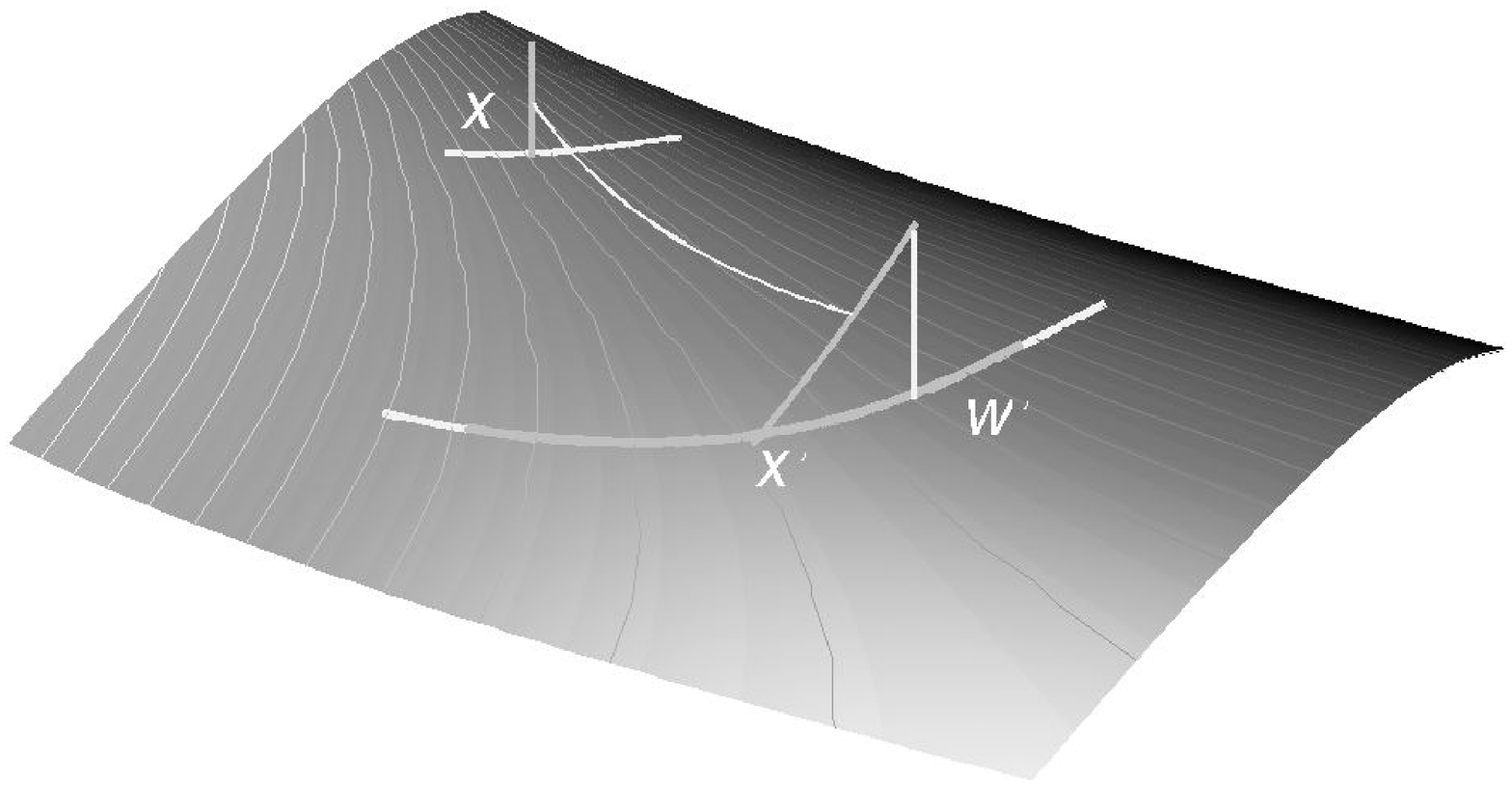}
\caption{Lemma \ref{holonomia no se tuerce sobre inestable}}
\label{figura inestable no se tuerce}
\end{figure}
\begin{lema}\label{holonomia no se tuerce sobre inestable}
There exists $k\in\Z^+$, not depending on $x$ nor on the choice of
the center curves, such that, in the hypotheses of Proposition
\ref{holonomia controlada}, the stable holonomy map $h^s$ from
$W^{cu}_{loc}(x)$ to $W^{cu}_{loc}(x')$ verifies
$$h^s(J^u_n(z))\subset J^{cu}_{n-k}(x')\qquad \forall n\geq k$$
for all $z\in B^c_n(x)$
\end{lema}
\bp Consider $x'\in \W^s_{loc}(x)$, and center curves
$W^c_{loc}(x)$, $W^c_{loc}(x)$ through $x, x'$ respectively,
contained in $W^{cs}_{loc}(x)$. Consider $y\in J^u_n(z)$,
with $z\in B^c_n(x)$, and let $y'=h^s(y)$, $z'=h^s(z)$. \par%
Let us find $k>0$ verifying:%
\begin{enumerate}%
\item $y'\in J^u_{n-k}(w')\subset J^{cu}_{n-k}(x')$ with %
\item $w'\in B^c_{n-k}(x')\subset W^c_{loc}(x')$%
\end{enumerate}
Since the point $f^n(y)$ is in $\W^u_{\nu_n(z)}(f^n(z))$, we have
$d(f^n(y),f^n(z))\leq \nu_n(z)$. Now, $y'\in \W^s_{loc}(y)$ and
$z'\in \W^s_{loc}(z)$, so:
$$d(f^n(y'), f^n(z'))\leq K\nu_n(y')\leq K^2\nu_n(z')$$
for a fixed constant $K>0$, not depending on $z$ (see Lemma \ref{holder cocycles} - Appendix 1).\par%
Let $w'\in \W^u_{loc}(y')\cap W^c_{loc}(x')$. From the fact that the
angle between the distributions is bounded from below, it follows by
projecting that
\begin{equation}\label{projection}
d(f^n(y'), f^n(w'))\leq C'\nu_n(y')\qquad \mbox{and}\qquad
d(f^n(w'), f^n(z'))\leq C'\nu_n(z')\end{equation} hence (1) follows
from the first inequality above by taking any $l_0>0$ verifying
$\nu_{-l_0}(y)>C'$ for all $y\in M$. Indeed,
$$d(f^{n-l}(y'),f^{n-l}(w'))\leq d(f^n(y'),f^n(w'))\leq C'\nu_n(y')\leq \nu_{n-l}(y')$$
for all $l\geq l_0$. Using Lemma \ref{holder cocycles} again, one
obtains $k>0$ such that $\nu_{-k}(y)>C$ for all $y\in M$, and so
$y'\in J^u_n(w')$.\par From the second inequality in
(\ref{projection}), and inequalities (\ref{adapted.metric.c}) and
(\ref{bunching condition}) in \S \ref{subsection definition of
juliennes} we derive
$$d(w', z')\leq C'\gamma_{-n}(z')\nu_n(z')\leq C'\s_n(z')\leq \s_{n-l}(z')$$
Now, previous lemma implies $z'\in B^c_{n-l}(x')$ for some
sufficiently large $l>0$, so using Lemma \ref{holder cocycles} again
and taking into account that $z'\in B^c_{n-l}(x')$, we find a
(uniform) $k>0$ so that $d(x',w')\leq \s_{n-k}(x')$ for all $n\geq
k$.\ep
\subsection{A characterization of Lebesgue density
points}\label{section cu density points}
In this paragraph, we shall see that the following three systems
are Vitali equivalent over essentially $u$-saturated sets:
\begin{enumerate}
\item $Q_n(x)=\bigcup_{y\in J^{sc}_n(x)}\W^u_{\s_n(y)}(y)$ where $J^{sc}_n(x)=\bigcup_{y\in B^c_n(x)}\W^s_{\s_{n}(y)}(y)$%
\item $J^{usc}_n(x)=\bigcup_{y\in J^{sc}_n(x)}J^u_n(y)$%
\item $J^{scu}_n(x)=\bigcup_{y\in J^{cu}_n(x)}\W^s_{\s_n(y)}(y)$
\end{enumerate}
The first system $Q_n(x)$ consists of ``cubic" balls, so it is not
difficult to see it is Vitali equivalent to Lebesgue. The second
system $J^{usc}_n(x)$ consists of dynamically defined local unstable
saturation of local center-stable leafs. Both systems are local
unstable saturations of the same center-stable leaf, and in both
cases the local unstable fibers are ``uniformly" sized, so over
essentially $u$-saturated sets, they have the same density points.
This is a consequence of absolute continuity of the unstable
foliation. Finally, the systems $J^{usc}_n(x)$ and $J^{scu}_n(x)$
are comparable, in the sense that they are nested, their volumes
preserving a controlled ratio. So the three systems are Vitali
equivalent over essentially $u$-saturated sets:
\begin{lema}\label{Q_n equivalente a cubos}
The system $\{Q_n(x)\}_{x\in M}$ is Vitali equivalent to Lebesgue.
\end{lema}
It follows from Proposition \ref{regular vitali equivalent} in
Appendix 2 and from the fact that the angle between the
distributions is bounded from below (note that all $x\in M$ verify
$\s_1^n\leq \s_n(x)\leq \s_2^n$ for some fixed
$\s_1,\s_2\in(0,1)$).\par We say that a measurable set $X$ is {\de
essentially $u$-saturated} if there exists a measurable
$u$-saturated set $X_u$ (an {\de essential $u$-saturate of $X$})
such that $m(X\triangle X_u)=0$.
\begin{prop}\label{USC_n equivalente a Q_n}
The system $\{J^{usc}_n(x)\}_{x\in M}$ is Vitali equivalent to
$\{Q_n(x)\}_{x\in M}$ over essentially $u$-saturated sets.
\end{prop}
\bp For measurable (small) sets $X$, let us denote by $m_u(X)$ and
$m_{sc}(X)$ the induced Riemannian volume of $X$ in $\W^u_{loc}$
and $W^{sc}_{loc}$ respectively (the choice of $W^{sc}_{loc}$ is
 fixed {\it a priori}). Since $\W^u$ is absolutely continuous, given
 any esentially $u$-saturated  $X$, and any essential $u$-saturate $X_u$ of $X$,
 we have
 \begin{enumerate}
 \item $m(X_u\cap Q_n(x))=\int_{X_u\cap J^{sc}_n(x)}m_u(\W^u_{\s_n(y)}(y))dm_{sc}(y)$%
 \item $m(X_u\cap J^{usc}_n(x))=\int_{X_u\cap J^{sc}_n(x)}m_u(J^u_n(y))dm_{sc}(y)$%
 \end{enumerate}
Observe that there exists a constant $D>1$ such that, for all $y\in
J^{sc}_n(x)$,
\begin{equation}\label{distorsion acotada de U_n}
\fra1D\leq \fra{m_u(J^u_n(y))}{m_u(J^u_n(x))}\leq D
\end{equation}
 (see lemma 4.1. of \cite{bw}).
Hence, we have,
$$\fra1{D^2}\fra{m_{sc}(X_u\cap J^{sc}_n(x))}{m_{sc}(J^{sc}_n(x))}\leq\fra{m(X_u\cap J^{usc}_n(x))}{m(J^{usc}_n(x))}\leq D^2\fra{m_{sc}(X_u\cap J^{sc}_n(x))}{m_{sc}(J^{sc}_n(x))}$$
And also,
$$\fra1{D^2}\fra{m_{sc}(X_u\cap J^{sc}_n(x))}{m_{sc}(J^{sc}_n(x))}\leq\fra{m(X_u\cap C_n(x))}{m(C_n(x))}\leq D^2\fra{m_{sc}(X_u\cap J^{sc}_n(x))}{m_{sc}(J^{sc}_n(x))}$$
So
$$\fra1{D^4}\fra{m(X\cap Q_n(x))}{m(Q_n(x))}\leq\fra{m(X\cap J^{usc}_n(x))}{m(J^{usc}_n(x))}\leq D^4\fra{m(X\cap Q_n(x))}{m(Q_n(x))}$$
The claim follows now from proposition \ref{regular vitali
equivalent}, part (3).\ep%
\begin{obs}
Observe that in the proof above we have used the same choice of
$B^c_n(x)$ for both $J^{usc}_n(x)$ and $Q_n(x)$; however, a fortiori
it follows that the choice of $B^c_n(x)$ is irrelevant.
\end{obs}
\begin{prop}
The system $\{J^{scu}_n(x)\}$ is Vitali equivalent to
$\{J^{usc}_n(x)\}$ over all measurable sets.
\end{prop}
\bp We shall find $\l\in{\mathbb Z}^+$ and $D>0$ such that
$$J^{scu}_{n+l}(x_0)\subset J^{usc}_n(x_0)\subset J^{scu}_{n-l}(x_0)\qquad \mbox{and}\qquad \fra{m(J^{usc}_{n+l}(x_0))}{m(J^{usc}_n(x_0))}\geq D$$
for all $x_0\in M$. The proof follows then from item (2) of
Proposition \ref{regular vitali equivalent}.\par Let us consider
$k_1>k$, where $k$ is the positive integer of Proposition
\ref{holonomia no se tuerce sobre inestable}, verifying $\min_{x\in
M}\s_{-k_1}(x)>C^2$ where $C$ is as in Lemma \ref{holder cocycles}.
If $z\in J^{usc}_n(x_0)$, then $z\in U_n(y)$, with $y\in
J^{sc}_n(x_0)$. By Lemma \ref{holonomia controlada sobre center} and
the choice of $k_1$, we have $y\in B^c_{n-k_1}(x)$, with $x\in
\W^s_{loc}(x_0)$. Applying Lemma \ref{holonomia no se tuerce sobre
inestable} to the holonomy map $h^s$ going from $J^{cu}_{n-k_1}(x)$
to $W^{cu}_{loc}(x_0)$, we have $ h^s(J^{cu}_{n-k_1}(x))\subset
J^{cu}_{n-2k_1}(x_0)$. Then, from the fact that the angles between
distributions is bounded from below, we
have that, for some $k_2>k_1$, $z\in J^{cu}_{n-k_1}(x)\subset J^{scu}_{n-k_2}(x_0)$\par%
The other inclusion is more simple, since, for $z\in
J^{scu}_n(x_0)$, we have $z\in \W^s_{\s_n(y)}(y)$ with $y\in
J^{cu}_n(x_0)$. But $W^{uc}_{loc}(z)\cap \W^s_{loc}(x_0)=\{x\}$, and
hence directly from lemma \ref{holonomia no se tuerce sobre
inestable} we have that $z$, belonging to $h^s(J^{cu}_n(x_0))$, is
contained in
$J^{cu}_{n-k_1}(x)$, hence $z\in J^{usc}_{n-k_1}(x_0)$.\par%
To finish the proof, let us see that
$m(J^{usc}_{n+l}(x))/m(J^{usc}_n(x))$ is bounded from below for all
$n>0$ and $x\in M$. Proceeding as in lemma \ref{USC_n equivalente a
Q_n}, we obtain that, there is a constant $c>0$ such that, for all
$x\in M$ and $n>0$
$$\fra1c\leq \fra{m(J^{usc}_n(x))}{m_u(J^u_n(x))m_s(\W^s_{\s_n(x)}(x))m_c(B^c_n(x))}\leq c$$
It is easy to see that
$m_s(\W^s_{\s_{n+l}(x)}(x))/m_s(\W^s_{\s_{n}(x)}(x))$ and
$m_c(B^c_{n+k}(x))/m_c(B^c_{n}(x))$ are uniformly bounded. Now, we
have
$$m_u(J^u_n(x))\leq K[{\rm Jac}(f^{-n})'(f^n(x))|_{E^u}]\lambda_n(x)$$
for some uniform $K>0$, so $m_u(J^u_{n+k}(x))/m_u(J^u_x(x))$ is
uniformly bounded too. For a detailed proof of this last
estimation see lemma 4.4 of \cite{bw}.\ep%
 {\em Proof of theorem \ref{saturacion de puntos de densidad}} %
 Let $X_s$ be an essential $s$-saturate of $X$. And assume $x$ is a $J^{scu}_n$- density point of
 $X$, hence of $X_s$.  Calling $m_s(A)$
the induced Riemannian volume of $A$ in $\W^s$, and $m_{cu}(A)$
the induced Riemannian volume of $A$ in some (fixed a priori)
$W^{cu}_{loc}$ we have, due to the fact that $X_s$ is
$s$-saturated:
$$\fra1K\leq \fra{m(X_s\cap J^{scu}_n(x))}{\s_n(x)m_{cu}(X_s\cap J^{cu}_n(x))}\leq K$$
Now, due to proposition \ref{holonomia controlada} we have
$$ m_{cu}(h^s(X_s\cap J^{cu}_{n+k}(x)))\leq m_{cu}(X_s\cap J^{cu}_{n}(h^s(x)))\leq m_{cu}(h^s(X_s\cap J^{cu}_{n-k}(x)))$$
The proof follows from the fact that
$$\fra1K\leq\fra{m_{cu}(h^s(X))}{m_{cu}(X)}\leq K$$
for some uniform $K>0$. \ep

\section*{Appendix 1}
\subsection*{H\"{o}lder cocycles and local leaves}
\setcounter{teo}{0}\setcounter{prop}{0}
\begin{lemaa}\label{holder cocycles}
For any H\"older continuous $\alpha:M\to\R^+$, there is a fixed
constant $C>1$ such that if $y\in \W^s_{loc}(x)\cup B^c_n(x)\cup
J^u_n(x)$, then $$\fra1C\leq \fra{\alpha_n(x)}{\alpha_n(y)}\leq
C\qquad \forall n\geq 0$$
\end{lemaa}
\bp See for instance \cite{bw2}\ep
\subsection*{Stable holonomy on center stable leaves}
The following is proved in particular in
\cite{bw2}:
\begin{propa} \cite{bw2} If $f:M\to M$ is a $C^{1+\alpha}$
partially hyperbolic with some center bunching condition
(trivially satisfied for one-dimensional center bundle), then
there exists $\beta>0$ such that the stable holonomy map between
center transversals is $C^{1+\beta}$
\end{propa}
We include a weaker version, for completeness, which is enough for
our purposes.
\begin{propa}\label{holonomia derivable}
There is a uniform Lipschitz constant $L>0$, such that for all
$x'\in \W^s_{loc}(x)$ and all central curves $W^c_{loc}(x)$,
$W^c_{loc}(x')$ contained in the same $W^{sc}_{loc}(x)$, the
stable holonomy map $h^s$ from $W^c_{loc}(x)$ to $W^c_{loc}(x')$
is $L$-Lipschitz when restricted to $W^{sc}(x)$. That is,
$$d(h^s(x),h^s(x'))\leq Ld(x,x')$$
\end{propa}%
We sketch the proof of this statement, the scheme of which may be
found in \cite{pughshub1972}. Take $W^{sc}(x)$, a $sc$ leaf
through $x$, $W^c_{loc}(x)\subset W^{sc}(x)$ a center curve. Take
also $W^{cu}_{loc}(y)$ a $cu$-leaf containing $W^c_{loc}(y)$ .

Now define $h^s:W^c_{loc}(x)\rightarrow W^{uc}_{loc}(y)$ in the
usual way (observe that $h^s(W^c_{loc}(x))\subset W^c_{loc}(y)$)

Take $S$ a smooth sub-bundle of $TM$ $C^0$ near to $E^s$. Fix
$\delta>0$ small and call $\mathcal S =\exp (\{s\in S_p\,
;||s||<\delta \})$. The map $p\mapsto \mathcal S$ is a smooth
pre-foliation.

Define the map $k_n:W^c_{loc}(x)\rightarrow W^{uc}_{loc}(y)$ in the
following way: take $z\in W^c_{loc}(x)$, $w= \mathcal S_{f^n(z)}\cap
f^n(W^{uc}_{loc}(h_s(z)))(\subset f^n(W^{uc}_{loc}(y)))$ and call
$k_n(z)=f^{-n}(w)$.\par
Observe that iteration for the past makes $f^{-n}(\mathcal S_p)$
converge uniformly on compact sets to $\W^s_{f^{-n}}(p)$ (and the
speed of convergence is independent of $p$). This observation easily
implies that $k_n$ uniformly converge to  $h^s$. Then, it is enough
to prove that $k_n$ are uniformly bounded for all $n$. The facts
that $f^n(W^c_{loc}(x))$ and $f^n(W^c_{loc}(y))$ are, increasing
$n$, as $C^1$ near as we want and that the angles between $S_p$ and
$E_{cu}$ are uniformly bounded from below,  give
 us that the map that sends $f^n(z)$ to $w$ is $C^1$-near to the
 inclusion. The uniform convergence of $f^{-k}(\mathcal S_p)$ to
 $\W^s_{f^{-k}}(p))$ again implies that, given $\varepsilon >0$
 there exists $n_0$ such that for all
 $n\geq n_0$ $d(f^{-k}(w),f^{-k}(f^n(z)))<\varepsilon$
$\forall 0\leq k\leq n$. Then, the Chain Rule and a typical
argument of distortion estimates of multiplicative H\"older
cocycles gives the desired bound for the derivative.

Finally observe that the Lipschitz constant only depends on the
stable distance between $x$ and $y$. This dependence appears in
the distortion estimates.

\section*{Appendix 2}
\subsection*{Vitali systems}
Let us briefly recall some known facts about density points. The
reader may see for instance \cite{shilov_gurevich}. We thank M.
Hirayama for pointing us a mistake in a previous statement of this
Proposition
\begin{propb}\label{regular vitali equivalent}
Each of the following are sufficient conditions for two systems
$\{B_n(x)\}_x$ and $\{C_n(x)\}_x$ to be Vitali equivalent over a
given $\sigma$-algebra ${\mathcal M}$:
\begin{enumerate}
\item There exist $k\in{\mathbb Z}^+$ and $D>0$ such that
$$B_{n+k}(x)\subset C_n(x)\subset B_{n-k}(x)\qquad\mbox{with}\qquad\fra{m(B_{n+k}(x))}{m(B_n(x))}\geq D\qquad \mbox{for all }x\in M $$
\item There exists $D>0$ such that
$$\fra1D\leq \fra{m(X\cap B_n(x))\,m(C_n(x))}{m(X\cap C_n(x))\,m(B_n(x))}\leq D \qquad \forall n\in{\mathbb Z}^+\quad \forall X\in {\mathcal M}$$
\end{enumerate}
\end{propb}


\end{document}